\documentclass{amsproc}
\usepackage[utf8]{inputenc}

\usepackage{amsmath,amssymb,amsthm,mathrsfs,mathtools}
\usepackage{hyperref}
\usepackage{xcolor}
\usepackage{pifont}

\usepackage{tikz}
\usepackage{tikz-cd}
\numberwithin{equation}{section}
\newtheorem{theorem}{Theorem}[section]

\theoremstyle{definition}
\newtheorem{definition}[theorem]{Definition}

\newtheorem{example}[theorem]{Example}

\newtheorem{algorithm}[theorem]{Algorithm}

\newcommand{\G}{\mathcal{G}}

\DeclareMathOperator{\sbridge}{sb}
\DeclareMathOperator{\lcm}{lcm}

\newcommand{\cmark}{\ding{51}}
\newcommand{\xmark}{\ding{55}}

\title{The MorseResolutions package for Macaulay2}

\author[T. Chau]{Trung Chau}
\address{Department of Mathematics, University of Utah, 155 South 1400 East, Salt Lake City, UT~84112, USA}
\email[Trung Chau]{trung.chau@utah.edu}

\author[S. Kara]{Selvi Kara}
\address{Science Research Initiative, University of Utah, 155 South 1400 East, Salt Lake City, UT~84112, USA}
\email[Selvi Kara]{selvi.kara@utah.edu}

\author[A. O'Keefe]{Augustine O'Keefe}
\address{Department of Mathematics and Statistics, Connecticut College\\
270 Mohegan Avenue Pkwy,
New London, CT 06320, USA}
\email[Augustine O'Keefe]{aokeefe@conncoll.edu}

\begin{document}

\begin{abstract}
Using discrete Morse theory, Batzies and Welker introduced Morse resolutions of monomial ideals. In this note, we present the {\it Macaulay2} package {\tt MorseResolutions} for working with two important classes of Morse resolutions: Lyubeznik and Barile-Macchia resolutions. This package also contains procedures to search for a minimal Barile-Macchia resolution of a given monomial ideal.
\end{abstract}


\maketitle

\section{Introduction}

A well-known approach for resolving a monomial ideal over a polynomial ring is to find its Taylor resolution \cite{Tay66} which are usually non-minimal. Several attempts have been made to find refinements of Taylor resolutions that are closer to the minimal resolutions. Among the most prominent refinements are Lyubeznik resolutions \cite{Ly88} and Scarf complexes \cite{BPS98}. Another type of refinement are \textit{Morse resolutions} which were introduced by Batzies and Welker in \cite{BW02}. These resolutions are induced from Taylor resolutions by utilizing \textit{homogeneous acyclic matchings}, a tool from discrete Morse theory.

In general, there are several homogeneous acyclic matchings for a given monomial ideal, and each such matching may produce a different Morse resolution. For instance, the homogeneous acyclic matchings described in \cite[Theorem 3.2]{BW02} yield a class of Morse resolutions which coincide with Lyubeznik resolutions.  In this note, we refer to these matchings as {\it Lyubeznik matchings}. Another class of Morse resolutions called \textit{Barile-Macchia resolutions} are induced by  a different class of  homogeneous acyclic matchings termed {\it Barile-Macchia matchings}, as described in \cite[Algorithm 2.9]{chau2022barile}. While Barile-Macchia resolutions may not always be minimal, they tend to be closer to the minimal resolution compared to the Taylor and Lyubeznik resolutions (in specific cases). Notably, they are minimal for many important classes of monomial ideals, as demonstrated in \cite{chau2022barile}.

Morse resolutions serve as powerful tools for describing (minimal) free resolutions of monomial ideals. To facilitate further exploration in this field, we have developed the {\tt MorseResolutions} package for {\it Macaulay2}, which enables researchers to experiment with Barile-Macchia and Lyubeznik resolutions. Specifically, our package allows users to generate Barile-Macchia matchings and Lyubeznik matchings for a given ideal, considering a specific total order on its minimal generating set. The package can be accessed at: \url{https://github.com/selvikara/morseResolutions}.

In this note, we provide a review of the necessary mathematical background and summarize the key features of our package, accompanied by illustrative examples.


\section{Mathematical background}
In this section, we summarize the relevant background for this note and package. We refer interested readers to \cite{chau2022barile} for more details on the theory.

Let $R=\Bbbk [x_1,\ldots, x_N]$ be the polynomial ring in $N$ variables over a field $\Bbbk$. For the remainder of this section, let $I \subset R$ be a monomial ideal with a minimal monomial generating set $\G(I)=\{m_1,\ldots, m_n\}$.  Impose a total order   $(>)$ on $\G(I)$ and treat a subset $\sigma$ of $\G(I)$ as an ordered set with respect to $(>)$. 

Throughout the note, let $[k]$ denote the set $\{1,2,\ldots, k\}$ for any integer $k$.

\subsection{Taylor resolutions}

Consider the full simplex on $n$ vertices such that the vertices are labeled elements of $\G(I)$  and each face is labeled by the lcm of its vertex labels. This labeled simplex, denoted by $\Delta_I$, is called the {\it Taylor complex} of $R/I$.  Let $\lcm(\sigma)$ denote the lcm of the elements in $\sigma \subseteq \G(I)$.

Recall that $\Delta_I$ induces a $\mathbb{Z}^N$-graded complex $\mathcal{F}$ where $\mathcal{F}_r$ is the free $R$-module with a basis indexed by all  subsets of cardinality $r$ of $\G(I)$ and the differentials $\partial_r\colon \mathcal{F}_r\to \mathcal{F}_{r-1}$ are  defined by
\[ \partial_r(\sigma) =\sum_{\substack{\sigma'\subseteq \sigma,\\|\sigma'|=r-1}} [\sigma:\sigma'] \frac{\lcm(\sigma)}{\lcm(\sigma')} \sigma'. \]
The complex $\mathcal{F}$ is called the {\it Taylor resolution} of $R/I$.

\subsection{Discrete Morse theory and Morse resolutions} 
Discrete Morse theory was developed by Forman in \cite{Fo94} as a combinatorial counterpart of Morse theory for manifolds, and was reformulated in terms of homogeneous acyclic matchings  by Chari in \cite{Cha00}. In \cite{BW02}, Batzies and Welker used Chari's reformulation to obtain a ``trimmed" cellullar resolution from a given cellular resolution induced by a regular CW-complex (for instance, the Taylor resolution) of a monomial ideal. The key idea of their method is to encode faces of the Taylor complex in a graph and then reduce this graph to obtain a smaller object that is homotopy equivalent to $\Delta_I$ by using homogeneous acyclic matchings.




 Consider the directed graph $G_I$ on the set of cells of the Taylor complex $\Delta_I$.  The vertices of $G_I$ are the cells of $\Delta_I$, namely subsets of $\G(I)$, and the directed edges of $G_I$ are given by
$$E_I= \{(\sigma,\sigma'): \sigma'\subseteq  \sigma, |\sigma'|=|\sigma|-1\}.$$
For any $A\subseteq E_I$, let $G_I^A$ be the graph obtained  by reversing the direction of  the directed edges in $A$, i.e., the directed graph with the edges
$$E(G_I^A)=(E_I \setminus A) \cup \{(\sigma', \sigma ) \mid (\sigma, \sigma') \in A\}.$$
We are now ready to define homogeneous acyclic matchings and critical cells.
\begin{definition}\label{def:acyclicmatch}
Let $A\subseteq E_I$ be a matching on $G_I$, i.e., no two edges in $A$ share a vertex. Then:
\begin{itemize}
    \item The matching $A$ is called {\bf acyclic} if the associated graph $G_I^A$ is acyclic, i.e., it does not contain any directed cycles.   
    \item The matching $A$ is called {\bf homogeneous} if $\lcm(\sigma)=\lcm(\sigma')$ for any directed edge $(\sigma, \sigma') \in A$.
\end{itemize}
If $A$ is a homogeneous acyclic matching on $G_I$, a cell $\sigma \in \G(I)$ is called $A$-{\bf critical} if $\sigma$ does not appear in any edge of $A$.
\end{definition}

In \cite{BW02}, Batzies and Welker showed that a homogeneous acyclic matching $A$ of $G_I$ produces a CW-complex which supports a free resolution of $R/I$. See \cite{BW02} and \cite{chau2022barile} for unexplained terminology and notations.

\begin{theorem}\cite{BW02} \label{thm:morseres}
If $A$ is a homogeneous acyclic matching on $G_I$, then there exists a CW-complex which supports a free resolution of $R/I$. The $i$-cells  of this CW-complex are in one-to-one correspondence with the $A$-critical cells of cardinality $i+1$ of $\Delta_I$.  We denote the cellular resolution supported on this CW-complex by $\mathcal{F}_A$ where $(\mathcal{F}_A)_i$ is the free $R$-module with a basis indexed by all critical cells   of cardinality $i$ and the differentials are the maps $\partial_i^A:(\mathcal{F}_A)_i\to (\mathcal{F}_A)_{i-1}$ defined by
     \[ \partial_i^A(\sigma) =\sum_{\substack{\sigma'\subseteq \sigma,\\|\sigma'|=i-1}} [\sigma:\sigma'] \sum_{\substack{\sigma'' \text{ critical,}\\ |\sigma''|=i-1 }} \sum_{\substack{\mathcal{P} \text{ gradient path}\\ \text{from } \sigma' \text{ to }\sigma''}} m(\mathcal{P}) \frac{\lcm(\sigma)}{\lcm(\sigma'')} \sigma''. \]
The resulting (cellular) free resolution $\mathcal{F}_A$ is called the {\bf Morse resolution} of $R/I$ associated to $A$.

The resolution $\mathcal{F}_A$ is minimal if for any two $A$-critical cells $\sigma, \sigma''$ of $\Delta_I$ with $|\sigma''|= |\sigma|-1$ such that there exists a gradient path from $\sigma'$ to $\sigma''$ in $G_I^A$ for some $\sigma'\subseteq \sigma$ with $|\sigma''|=|\sigma'|$, we have $\lcm(\sigma)\neq \lcm(\sigma'')$.
\end{theorem}

\subsection{Lyubeznik resolution} In \cite{BW02}, Batzies and Welker showed that Lyubeznik's resolution from \cite{Ly88} is an example of a Morse resolution and they provided a recipe for the corresponding homogeneous acyclic matching. We recall the construction below.  In the {\tt MorseResolutions} package, we provide a function that produces this Lyubeznik matching of $R/I$ with respect to a given total order.

\begin{theorem}\cite[Theorem 3.2]{BW02}\label{def:Lyuusingmatchings}
   For any subset $\sigma=\{m_1, \ldots ,m_q\}$ of $\G(I)$ where $m_1> \cdots > m_q$, we define
    \[
    	v_L(\sigma)\coloneqq  \sup \big\{ k\in \mathbb{N} : \exists m \in \G(I) \text{ such that } m_k > m 
    	\text{ for some } k\in [q] \text{ and } m\mid \lcm(m_1, \ldots, m_k) \big\}.
  \]
    Set $v_L(\sigma)= -\infty$ if no such $m$ exists. If $v_L(\sigma)\neq -\infty$, define
     \[
    m_L(\sigma)\coloneqq \min_{>} \{m\in \G(I): m\mid \lcm(m_1,\dots, m_{v_L(\sigma)})  \}.
    \]
    For each $p \in R$, set
    \[
    A_p\coloneqq \{(\sigma\cup m_L(\sigma), \sigma \setminus m_L(\sigma)) :  \lcm (\sigma)=p \text{ and }v_L(\sigma)\neq -\infty \}.
    \]
    Then  $\displaystyle A=\cup_{p\in R }A_p$ is a homogeneous acyclic matching. The graded free resolution $\mathcal{F}_A$ induced by $A$ is the {\bf Lyubeznik resolution} of $R/I$ with respect to $(>)$.
\end{theorem}

We call the final matching $A$ in this theorem the {\bf Lyubeznik matching} of $R/I$ with respect to $(>)$. Recall the following example from \cite[Example 5.6]{chau2022barile} to see how to find the $v_L(\sigma)$ and $m_L(\sigma)$ for a given  $\sigma$ of $\G(I)$:
\begin{example}\label{ex:BWmatching}
      Consider the monomial ideal $I=(m_1,m_2,m_3,m_4,m_5,m_6) \subseteq \Bbbk [x_1, \ldots, x_8]$ where
\[
		m_1=x_1x_2x_3x_4, ~~m_2=x_2x_3x_5x_6, ~~m_3=x_1x_2x_5, ~~m_4=x_1x_2x_7, ~~m_5=x_2x_3x_8, ~~m_6=x_7x_8
\]
with the total order $ m_1 > m_2 > m_3 > m_4 > m_5 > m_6$  on $\G(I)$. 
\begin{itemize}
    \item For $\sigma_1=(m_1,m_4,m_5)$, we have $v_L(\sigma_1)=3$ and $m_L(\sigma_1)=m_6$ since  $m_6 | \lcm(\sigma_1)$ and $\lcm(m_1,m_4)$ is not divisible by $m_5$ or $m_6$.
    \item For  $\sigma_2=(m_1,m_2,m_3)$, we have $v_L(\sigma_2)=2$ and $m_L(\sigma_2)=m_3$ since $m_3 |  \lcm(m_1,m_2)$.
    \item For $\sigma_3=(m_2,m_3,m_4)$, we have $v_L(\sigma_3)=- \infty$ since  $ \lcm(m_2,m_3)$ is not divisible by $m_4, m_5$ or $m_6$ and $\lcm(\sigma_3)$ is not divisible by $m_5$ or $m_6$.
\end{itemize}
Let $A$ be the Lyubeznik matching of $R/I$ with respect to this total order. Then,  $(\sigma_1 \cup \{m_6\}, \sigma_1)$ and $(\sigma_2, \sigma_2\setminus \{m_3\})$ are directed edges in $A$ while $A$ does not have any edge involving $\sigma_3$.
\end{example}

\subsection{Barile-Macchia resolution} The first two authors developed an algorithm to produce a different homogeneous acyclic matching called the {\bf Barile-Macchia matching} of $R/I$ in \cite{chau2022barile}.  Morse resolutions induced by Barile-Macchia matchings are called {\bf Barile-Macchia resolutions}. In the same work, they identified several classes of ideals whose Barile-Macchia resolutions are minimal. In {\tt MorseResolutions} package, we provide a function that produces the Barile-Macchia matching of $R/I$ with respect to a given total order.

\begin{definition}
Let $\sigma$ be a subset of $\G(I)$. 
\begin{itemize}
    \item A monomial $m\in \G(I)$ is called a \textit{bridge} of $\sigma$ if $m\in \sigma$ and $\lcm (\sigma \setminus \{m\})=\lcm (\sigma)$.
    \item  The \emph{smallest bridge} of $\sigma \subseteq \G(I)$ with respect to $(>)$, denoted by $\sbridge_{>}(\sigma)$, is the smallest bridge of $\sigma$ in this total order.  If $\sigma$ has no bridges, then $\sbridge_{>}(\sigma) = \emptyset$.
\end{itemize}
\end{definition}

For simplicity, denote $S \setminus \{ s\}$ and $S\cup \{s\}$ for a set $S$ and $s\in S$ by $S \setminus s$ and $S\cup s$, respectively. Below is the Barile-Macchia algorithm from \cite[Algorithm 2.9]{chau2022barile}.
 
\begin{algorithm}[Barile-Macchia Algorithm]\label{alg}
    {\sf Let $A=\emptyset$. Set 
    $$\Omega=\{\text{all  subsets of } \G(I) \text{ with cardinality at least } 3\}.$$
    \begin{enumerate}
        \item Pick a subset $\sigma$ of maximal cardinality in $\Omega$. 
        \item  Set
        \[ \Omega \coloneqq \Omega \setminus \{\sigma, \sigma \setminus \sbridge(\sigma)\}. \]
        If  $\sbridge(\sigma)\neq \emptyset$, add the directed edge $(\sigma , \sigma \setminus \sbridge(\sigma))$ to $A$. If $\Omega\neq \emptyset$, return to step (1). 
        \item Whenever there exist distinct directed edges $(\sigma, \sigma \setminus \sbridge(\sigma))$ and $(\sigma', \sigma' \setminus \sbridge(\sigma'))$ in $A$ such that 
            $$\sigma \setminus \sbridge(\sigma) = \sigma' \setminus \sbridge(\sigma'),$$
            then 
            \begin{itemize}
                \item if $ \sbridge(\sigma') > \sbridge(\sigma)$, remove $(\sigma', \sigma' \setminus \sbridge(\sigma'))$ from $A$,
                \item  otherwise, remove $(\sigma, \sigma \setminus \sbridge(\sigma))$ from $A$.
            \end{itemize}
    \end{enumerate}}
\end{algorithm}

\begin{definition}
  The set $A$ produced by the Barile-Macchia algorithm is called the {\bf Barile-Macchia matching} of $R/I$ with respect to $(>)$.  The set of edges that appear in $A$ before Step (3) of the algorithm is called the {\bf possible edges} of $R/I$. 
\end{definition}
   Barile-Macchia matchings are indeed homogeneous acyclic matchings as shown in \cite[Theorem 2.11]{chau2022barile}.  Recall the following example from \cite[Example 2.13]{chau2022barile}.
\begin{example}\label{ex:1}
    Consider the ideal $I=(wz,wx,xy,yz) \subseteq R=\Bbbk[w,x,y,z]$ with the total order $wz> wx> xy> yz$ on $\G(I)$. The Barile-Macchia matching of  $R/I$ with respect to this total order is:  
    $$A=\{(\{wz,wx,xy,yz\},\{wz,wx,xy\}), (\{wz,xy,yz\}, \{wz,xy\}), (\{wx,xy,yz\},\{wx,yz\})\}.$$ 
   The possible edges of $A$ are 
   $$A \cup \{(\{wz,wx,yz\}, \{wx,yz\})\}$$
   
  \noindent because the edge $(\{wz,wx,yz\}, \{wx,yz\})$ is removed from $A$ in Step (3) of the algorithm. We will display the same output by using our functions from the {\tt MorseResolutions} package.
\end{example}

A final concept that is needed in this section is the notion of bridge-friendliness from \cite[Definition 2.27]{chau2022barile}. As it is shown in \cite[Theorem 2.29]{chau2022barile},  if an ideal is bridge-friendly, then its corresponding Barile-Macchia resolution is minimal. So, this property allows us check the minimality of Barile-Macchia resolutions. In this package, we provide a function to identify whether a given ideal is bridge-friendly while detecting all total orders where this property holds. 
\begin{definition}
   A monomial ideal $I$ is called {\bf bridge-friendly} with respect to a total order $(>)$ on $\G(I)$ if  its  Barile-Macchia matching with respect to $(>)$ is the same as its set of  possible edges.
\end{definition}

The ideal in Example \ref{ex:1} is not bridge-friendly with respect to the given order. More generally, we will show that this ideal is not bridge-friendly with respect to any total order by using our package in the next section.

                        	
\section{The package and examples}

The {\it Macaulay2} package {\tt MorseResolutions} was
created as a tool to experiment with Barile-Macchia resolutions from \cite{chau2022barile}, Lyubeznik resolutions from \cite{BW02}, and the notion of bridge-friendliness from \cite{chau2022barile}. In this section, we highlight some of the key features of this package.

Our package has two main parts where the first part focuses on functions related to Barile-Macchia resolutions and the second part involves analogous functions for Lyubeznik resolutions. This package also contains an application of these ideas, trimmed Lyubeznik resolutions, in Subsection \ref{sub:4.2}.

\subsection{Barile-Macchia matchings and Barile-Macchia resolutions}

One of the main functions of the package is the {\tt bMMatching} function which produces a Barile-Macchia matching of a given ideal with respect to a fixed total order $(>)$ based on Algorithm \ref{alg}. As seen from Theorem \ref{thm:morseres}, finding the Barile-Macchia matching is the key step in obtaining the corresponding Barile-Macchia resolution of the ideal.

In order to find a Barile-Macchia (or Lyubeznik) matching of $R/I$ for a given ideal $I$, the user must specify a total order on $\G(I)$. In our package, we record a total order $(>)$ on $\G(I)$  as an ordered sequence $S$ such that  elements of $S$ are ordered from the smallest to the largest with respect to $(>)$. 

\begin{example}\label{ex:package1}
    Consider the ideal $I=\langle wx,xy,yz,wz \rangle \subseteq \Bbbk [w,x,y,z]$ from Example \ref{ex:1} with the total order $wz>wx>xy>yz$. The ordered set for this total order is $S=(yz,xy,wx,wz)$.
\end{example}

For a given total order $S$ on $\G(I)$, the function {\tt bMMatching} returns a list which consists of all the edges of the Barile-Macchia matching of $R/I$ with respect to $S$. As an illustration, we consider the ideal $I$ from Example \ref{ex:package1} with the total order $S=(yz,xy,wx,wz)$. Below is the Barile-Macchia matching of $R/I$ with respect to $S$.

\begin{verbatim}
i1 : R=QQ[w,x,y,z];
i2 : S = (y*z, x*y, w*x, w*z);
i3 : bMMatching(S)
o3 = {{({y*z, x*y, w*x}, {y*z, w*x}), ({y*z, x*y, w*z}, {x*y, w*z}), ({y*z, x*y, w*x,
      --------------------------------------------------------------------------------
      w*z}, {x*y, w*x, w*z})}}
\end{verbatim}

Recall that the collection of edges at the end of Step (2) of the algorithm is called the possible edges of $I$ with respect to $S$. In Step (3) of the algorithm, some of the edges from this colllection are removed to guarantee the output of the algorithm is a matching. So, the first step in obtaining a Barile-Macchia matching of $R/I$ with respect to $S$ is to identify its possible edges.

For a given ideal $I$, the function {\tt possibleEdges} returns a list of possible edges of $R/I$ with respect to $S$. Possible edges of $R/I$ with respect to $S$ from Example \ref{ex:package1} is as follows:

\begin{verbatim}
i4 : possibleEdges(S)
o4 = {({y*z, x*y, w*x, w*z}, {x*y, w*x, w*z}), ({y*z, x*y, w*x}, {y*z, w*x}), 
      --------------------------------------------------------------------------------
    ({y*z, x*y, w*z}, {x*y,  w*z}), ({y*z, w*x, w*z}, {y*z, w*x})}
\end{verbatim}

The {\tt possibleEdges} calls  {\tt possibleEdgesWithPositions} which returns a list of triples $(i, \sigma, \sigma \setminus \sbridge(\sigma))$, where $i$ is the position of the smallest bridge of  $\sigma \subseteq \G(I)$ with respect to a total order $S$ on $\G(I)$. In addition, the {\tt bMMatching} builds on {\tt possibleEdgesWithPositions}. This is because one must be able to compare the positions of the smallest bridges of two cells $\sigma$ and $\sigma'$ where $\sigma \setminus \sbridge(\sigma) =\sigma' \setminus \sbridge(\sigma')$ in  Step (3) of the algorithm.   For this reason,  we record the positions of the smallest bridges of cells in the output of this function.
\begin{verbatim}
i5 : possibleEdgesWithPositions(S)
o5 = {(0, {y*z, x*y, w*x, w*z}, {x*y, w*x, w*z}), (1, {y*z, x*y, w*x}, {y*z, w*x}), 
      --------------------------------------------------------------------------------
    (0, {y*z, x*y, w*z}, {x*y,  w*z}), (3, {y*z, w*x, w*z}, {y*z, w*x})}
\end{verbatim}

If the user wants to check whether the possible edges and the Barile-Macchia matching of $R/I$ with respect to $S$ overlap, we have built the function {\tt isBridgeFriendly} to test it. This function calls {\tt possibleEdges} and {\tt bMMatching} and then returns a Boolean. As we can see from the outputs of these two functions, our running example is not bridge-friendly.

\begin{verbatim}
i6 :  isBridgeFriendly(S)
o6 =  false
\end{verbatim}

We can also search over all possible total
orders on $\G(I)$ to determine the set of total orders such that $I$ is bridge-friendly. Specifically, the command {\tt bridgeFriendlyList}
returns  a list of pairs, where the first element in the pair is the total order for which the
ideal $I$ is bridge-friendly and the second element is the Barile-Macchia matching of $R/I$ with respect to that order. Note that running this command is computationally expensive as it  checks $N!$ different total orders in $R=\Bbbk [x_1,\ldots, x_N]$.

As remarked earlier, the ideal in our running example is not bridge-friendly with respect to any of the total orders on $\G(I)$. 

\begin{verbatim}
i7 : bridgeFriendlyList(I)
o7 = {}
\end{verbatim}

To illustrate a different output of the {\tt bridgeFriendlyList} function, we consider the following ideal $I$ which is bridge-friendly with respect to all total orders on $\G(I)$: 

\begin{verbatim}
i8 : I=ideal(x*y,y*z,x*z);
i9 : bridgeFriendlyList(I)
o9 = {{{x*y, y*z, x*z}, {({x*y, y*z, x*z}, {y*z, x*z})}}, 
------------------------------------------------------------------
      {{x*y, x*z, y*z}, {({x*y, x*z, y*z}, {x*z, y*z})}},
------------------------------------------------------------------
      {{y*z, x*y, x*z}, {({y*z, x*y, x*z}, {x*y, x*z})}}, 
------------------------------------------------------------------
      {{y*z, x*z, x*y}, {({y*z, x*z, x*y}, {x*z, x*y})}},
------------------------------------------------------------------
      {{x*z, x*y, y*z}, {({x*z, x*y, y*z}, {x*y, y*z})}}, 
------------------------------------------------------------------
      {{x*z, y*z, x*y}, {({x*z, y*z, x*y}, {x*y, y*z})}}}
\end{verbatim}

For a given ideal $I$ with respect to a total order $S$, the function {\tt criticalBMCells}  returns a list of all $A$-critical cells of $I$ where $A$ is the Barile-Macchia matching of $R/I$ with respect to $S$. This
function calls {\tt bMMatching} and the list goes through all Taylor cells of $I$ from larger cardinality towards smaller ones. Consider the ideal $I$ with its total order $S$ from Example \ref{ex:package1}.  In {\tt o9}, the first element in the list is the empty set. This means the Taylor cell of maximum cardinality appeared in the Barile-Macchia matching.

\begin{verbatim}
i10 : criticalBMCells(S)
o10 = {{}, {{y*z, w*x, w*z}}, {{y*z, x*y}, {x*y, w*x}, {y*z, w*z}, {w*x, w*z}}, 
      --------------------------------------------------------------------------------
{{y*z}, {x*y}, {w*x}, {w*z}}}
\end{verbatim}

Note that none of the critical cells  in {\tt o9} have the same $\lcm$. Thus, the corresponding Barile-Macchia resolution is minimal by Theorem \ref{thm:morseres}. 

One can deduce the ranks of the free $R$-modules in the Barile-Macchia resolution of $R/I$ with respect to $S$ from the critical cells of $I$ by Theorem \ref{thm:morseres}. The $i^{th}$ element in this list is the rank of the free $R$-module at homological degree $i$. Ranks are read from left to right following the tradition of {\it Macaulay2}. Consider  the ideal $I$ with the total order  $S$ from Example \ref{ex:package1}. Notice that {\tt bMRanks(S)} coincides with the total Betti numbers of $R/I$. 
\begin{verbatim}
i11 : bMRanks(S)
o11 = {1, 4, 4, 1, 0}
i12 :  res I
       1      4      4      1
o12 = R  <-- R  <-- R  <-- R  <-- 0 
      0      1      2      3      4
\end{verbatim}
\subsection{Lyubeznik matchings and resolutions} Similar to the functions for Barile-Macchia matchings and  resolutions, this package contains  {\tt lyubeznikMatching}, {\tt criticalLyubeznikCells(S)}, and  {\tt lyubeznikRanks(S)} functions to produce Lyubeznik matching, critical Lyubeznik cells, and the ranks of the free $R$-modules in the Lyubeznik resolution of a given ideal with respect to a fixed total order $(>)$.  For our running example $I$ and $S$, we display the outputs of these functions below. 
\begin{verbatim}
i13 : lyubeznikMatching(S)
o13 = {({y*z, x*y, w*x, w*z}, {x*y, w*x, w*z}), ({y*z, x*y, w*z}, {x*y, w*z})}
      --------------------------------------------------------------------------------
{{y*z}, {x*y}, {w*x}, {w*z}}}

i14 : criticalLyubeznikCells(S)
o14 = {{}, {{y*z, x*y, w*x}, {y*z, w*x, w*z}}, {{y*z, x*y}, {y*z, w*x}, {x*y, w*x}
 --------------------------------------------------------------------------------
 {y*z, w*z}, {w*x, w*z}}, {{y*z}, {x*y}, {w*x}, {w*z}}}, {{y*z}, {x*y},{w*x}, {w*z}}}
 
i15 : lyubeznikRanks(S)
o15 = {1, 4, 5, 2, 0}
\end{verbatim}

As it can be seen from {\tt o13} and Theorem \ref{thm:morseres}, the Lyubeznik resolution of $R/I$ with respect to $S$ is not minimal because the cells $\sigma=\{yz, xy, wx\}$ and $\sigma'=\{yz, wx\}$ are critical even though  $\sigma' \subset \sigma$ and $\lcm (\sigma)=\lcm (\sigma')$. We discuss how to trim this Lyubeznik resolution to obtain a resolution closer to the minimal one in Subsection \ref{thm:morseres}.

\section{Application}

In this section, we provide two applications of the {\tt MorseResolutions} package.

\subsection{Edge ideals of cycles} Using our package {\tt MorseResolutions} we found minimal examples of graphs whose edge ideals are not bridge-minimal. A monomial is called \textit{bridge-minimal} if there exists a total order $S$ such that its corresponding Barile-Macchia resolution is minimal. The bridge-minimality concept was introduced in \cite{chau2022barile}. Among cycles, the edge ideal of a $9$-cycle is the first cycle that is neither bridge-friendly nor bridge-minimal.

Table \ref{table1}, first appeared in \cite{chau2022barile}, summarizes the results
of our computations for bridge-friendliness and bridge-minimality of edge ideals of cycles up to 10 vertices. The first column of Table \ref{table1} is the number of vertices, the second column displays whether the edge ideal of the cycle is bridge-friendly with respect to some total order while the last column records whether any of the Barile-Macchia resolutions of the edge ideal of a cycle is minimal.

\begin{table}[h]
\centering
    \begin{tabular}{|c|c|c|} 
        \hline
        \hline
        Number of vertices & Bridge-friendliness & minimal Barile-Macchia resolution\\ \hline
        3,5,6 & \cmark      & \cmark         \\ \hline
        4,7,8,10 & \xmark      & \cmark         \\ \hline
        9 & \xmark      & \xmark         \\ \hline
        \hline
    \end{tabular}
    \vspace{.2cm}
    \caption{Bridge-friendliness and minimality of Barile-Macchia resolutions of edge ideals of cycles}
    \label{table1}
\end{table}

 This data was computed also using the Macaulay2 package EdgeIdeals \cite{francisco2009edgeideals}. We used the following code to identify edge ideals of cycles that are not bridge-friendly:

    \begin{verbatim}
loadPackage "EdgeIdeals"    
for i from 3 to 9 do (
R = QQ[x_1..x_i];
Ci = cycle(R,i);
I = edgeIdeal Ci;
if bridgeFriendlyList(I)=={} then print (i);
);
    \end{verbatim}
    Since the Betti numbers of edge ideals of cycles do not depend on characteristic of the base field $\Bbbk$ \cite{Jac04}, we can check the bridge-minimality of these ideals by comparing their Betti numbers with the ranks of its Barile-Macchia resolutions. The following code is used to check the bridge-minimality of edge ideal of a 9-cycle, which can be modified to check for other cycles:
   \begin{verbatim}
R=QQ[x1,x2,x3,x4,x5,x6,x7,x8,x9]    
I=ideal(x1*x2,x2*x3,x3*x4,x4*x5,x5*x6,x6*x7,x7*x8,x8*x9,x1*x9)
orderings = permutations flatten entries gens I;
st = false;
ij = 0;
while (st==false) do(
    if (ij < #orderings) then (
        ranks=bMRanks(toSequence orderings_ij);
        if (ranks=={1, 9, 27, 39, 27 ,9 , 2 , 0, 0, 0}) then(
            st=true; 
            print "This ideal is bridge-minimal.";
            );
            )
        else (
            st=true;
            print "This ideal is not bridge-minimal.";
            );
       ij = ij+1;     
    );
    \end{verbatim}
\subsection{Trimming Lyubeznik resolutions}\label{sub:4.2} In \cite{chau2022barile}, the Barile-Macchia matching algorithm is applied to the Taylor complex of a given monomial ideal with respect to a fixed total order. In this set-up, the Taylor complex can be replaced with any regular CW-complex (see \cite{BW02}). Since Lyubeznik resolutions are simplicial by \cite{Mer09} and any simplicial complex is a regular CW-complex, one can apply the Barile-Macchia algorithm to the simplicial complex that supports a Lyubeznik resolution of a given monomial ideal. 

In  {\tt MorseResolutions} package, we introduced the {\tt trimmedMatching} function, a modified version of the {\tt bMMatching} function. Instead of Taylor cells, this new function considers the critical Lyubeznik cells of a monomial ideal with respect to a total order $S_1$ and applies the Barile-Macchia algorithm  to these cells with respect to $S_2$. Consider our running example and let $S_1=S_2=S=(yz,xy,wx,wz)$. The output of this function is the collection of edges in $G_I$ with  ``unmatched" vertices in the critical Lyubeznik cells. 
\begin{verbatim}
i15 : trimmedMatching(S,S)
o15 =  {({y*z, x*y, w*x}, {y*z, w*x})}
\end{verbatim}
In addition, the package contains {\tt criticalTrimmedCells} and  {\tt trimmedRanks} functions to print out the critical cells of the trimmed Lyubeznik complex and the ranks of the trimmed Lyubeznik resolution. For our running example, the trimmed Lyubeznik resolution of $R/I$ is minimal.
\begin{verbatim}
i16 : criticalTrimmedCells(S,S)
o16 =  {{}, {{y*z, w*x, w*z}}, {{y*z, x*y}, {x*y, w*x}, {y*z, w*z}, {w*x, w*z}}
      --------------------------------------------------------------------------------
{{y*z}, {x*y}, {w*x}, {w*z}}}

i17 : trimmedRanks(S,S)
o17 =  {1, 4, 4, 1, 0}
\end{verbatim}
In fact, for any total orders $S_1$ and $S_2$, the corresponding trimmed Lyubeznik resolution of $R/I$ is minimal.

\noindent
{\bf Acknowledgments.} 
The authors thank Jeremy Dewar for helpful discussions. Chau was supported by NSF grants DMS 1801285, 2101671, and 2001368.   

\bibliographystyle{abbrv}
\bibliography{refs}

\end{document}